\documentclass[review,12pt]{elsarticle}

\usepackage{amsmath}
\usepackage{amssymb}
\usepackage{amsthm}
\usepackage{hyperref}
\hypersetup{
	colorlinks=true,
	citecolor=blue,
	linkcolor=blue
}
\usepackage[T1]{fontenc}
\usepackage[utf8]{inputenc}

\newtheorem{theorem}{Theorem}[section]
\newtheorem{lemma}[theorem]{Lemma}
\newtheorem{corollary}[theorem]{Corollary}
\newtheorem{remark}[theorem]{Remark}
\newtheorem{definition}[theorem]{Definition}
\newtheorem{openproblem}[theorem]{Open Problem}

\journal{J. Math. Anal. Appl.}
\begin{document}
\begin{frontmatter}
\title{Completely monotone functions in general and some applications}
\author[inst1]{V.E.S. Szabó}
\ead{bmesszabo@gmail.com}
\affiliation[inst1]{organization={Department of 
Analysis and Operations Research, 
Institute of Mathematics, \\
Budapest University of Technology and 
Economics}, 
            addressline={\\ Műegyetem rkp. 3.},
            city={Budapest},
             postcode={H-1111},
             country={Hungary}}         
\begin{abstract}
We simplify the proof of some widely used theoretical theorems, extending their 
applicability, while correcting some erroneous results. We also generalise key
results and present new results that contribute to the development of the
theory. Furthermore, we use the results obtained to investigate the
monotonicity properties of some specific functions related to the Gamma function. 
Finally, we formulate an open problem related to the psi function.
\end{abstract}            
\end{frontmatter}

\section{Introduction}
\label{sec:Intro}

The subject of this paper is completely monotone and related functions. 
First of all, we recall some definitions.

\begin{definition}
A function $f:I\to \mathbf{R}$, where $I\subseteq\mathbf{R}$ is an 
interval, is called completely monotone if $f$ has derivatives of
all orders and satisfies
\begin{equation*}
	(-1)^n f^{(n)}(x) \geq 0\quad (n=0,1,2,\ldots;\, x\in I).
\end{equation*}
\end{definition}

\begin{definition}
A positive function $f:I\to \mathbf{R}$, where $I\subseteq\mathbf{R}$ is 
an interval, is called logarithmically completely monotone if $f$ has 
derivatives of all orders and satisfies
\begin{equation*}
	(-1)^n (\log\circ f)^{(n)}(x) \geq 0\quad (n=1,2,\ldots;\, x\in I).
\end{equation*}
\end{definition}

\begin{definition}
A function $f:I\to \mathbf{R}$, where $I\subseteq\mathbf{R}$ is an 
interval, is called absolutely monotone if $f$ has derivatives of
all orders and satisfies
\begin{equation*}
	f^{(n)}(x) \geq 0\quad (n=0,1,2,\ldots;\, x\in I).
\end{equation*}
\end{definition}

\begin{definition}
A function $f:I\to \mathbf{R}$, where $I\subseteq\mathbf{R}$ is an
interval, is called a Bernstein function if $f\geq 0$, and $f'$ is 
completely monotone.
\end{definition}

It is obvious that a function $f(x)$ is completely monotone in $(a,b)$ if and
only if $g(x):=f(-x)$ is absolutely monotone in $(-b,-a)$. 
Also, 
a positive function $f$ is logarithmically completely monotone 
if and only if $(-\log\circ f)'$ is completely monotone. 
We note here that if $I=(0,\infty)$, $f\geq 0$, $f'\leq 0$, 
and $(-1)^n f^{(n)}\geq 0$ for infinitely many of $n\in\mathbf{N}$,  
then $f$ is completely monotone; see \cite{Gneiting1998}.

Some authors introduce the notion of strictly completely monotone, etc., 
if we have strict inequality in the definitions, however if $I=(a,\infty)$, 
and $f:(a,\infty)\to\mathbf{R}$ is completely monotone where 
$a$ is finite or $a=-\infty$, and $f$ is not 
a constant function, then
\begin{equation*}
	(-1)^n f^{(n)}(x) > 0\quad (n=0,1,2,\ldots;\, x\in (a,\infty)),
\end{equation*}
see \cite{Dubourdieu1939}, and see also \cite{LorchSzego1963} Section 8.

The main motivation of this paper is to discuss the wrong proof (the "only if" 
part) of 
\cite{Chen2007} Theorem 1 
(Theorem \ref{thm:psi_Chen2007_improve_comp_mon} $(iii_2)$ and Remark \ref{rem:alpha_less_2}), 
to correct the wrong proof of 
\cite{Chen2007} Lemma 1 (Lemma \ref{lem:phi_x_per_x}), 
to give simpler proof of \cite{QiChen2004} Theorem 2 (Theorem 
\ref{thm:log_gamma_vogt_general}) and \cite{Chen2007} Theorem 2 
(Theorem \ref{thm:wrong_mon}), and to generalize them, and   
to generalize \cite{Chen2007} Theorem 3 (Theorem \ref{thm:frac_gamma_power}),  
each of them was published in this journal.

Furthermore, we give simpler proof of known theorems 
(Theorem \ref{thm:product_of_comp_mon}, Theorem 
\ref{thm:composition_comp_mon_1}, Corollary 
\ref{cor:log_f_diff_comp_mon}, Corollary \ref{cor:logcomp_comp}, 
Lemma \ref{lemma:1perx_comp_mon}), 
extend the 
validity of them 
(Theorem \ref{thm:analytic_abs}, Corollary \ref{cor:analytic_comp}), 
generalize them
(Theorem \ref{thm:fractionpower_log_comp_mon}, 
Theorem \ref{thm:psi_Chen2007_improve_comp_mon} $(iii_1)$), 
and prove some new results
(Theorem \ref{thm:composition_comp_mon_2}, 
Theorem \ref{thm:f_powerto_g_log_comp}, 
Theorem \ref{thm:linfraction_log_comp_mon}, 
Theorem \ref{thm:fraction_of_Gamma_log_comp}). 
More details are given below.

In \cite{Widder1946} it was proved that if a function is absolutely monotone 
in $a\leq x<b$ then it can be extended analytically into the complex plane 
and the function will be analytic in the circle $|z-a|<b-a$. In fact, the
weaker assumption $a<x<b$ is enough in Theorem \ref{thm:analytic_abs}, and
we obtain the corresponding statement for the completely monotonic function
in Corollary \ref{cor:analytic_comp}. 

Proving that the product of completely monotone functions is completely 
monotone, the Leibniz formula for the $n$-th derivative of the product of 
functions can be applied. Here we give a simpler, 
very short and direct proof in 
Theorem \ref{thm:product_of_comp_mon}.

A nice and useful corollary of Theorem \ref{thm:product_of_comp_mon} is Corollary 
\ref{cor:f_power_g_logcomp_mon}, we will use it in the proof of 
Theorem \ref{thm:f_powerto_g_log_comp}.

The next theorem, Theorem \ref{thm:composition_comp_mon_1},  
 can be found in \cite{LorchNewman1983} saying that 
 Leibniz's formula for the derivative of product implies the statement, 
but omitting the details. An alternative proof using a formula for the 
$n$-th derivative of a composite function can be found in 
\cite{MillerSamko2001} Theorem 2. Here we show that using Leibniz's formula 
the statement follows. An easy consequence is Corollary 
\ref{cor:f_x_alpha}. The other consequence is Corollary 
\ref{cor:log_f_diff_comp_mon} that is not new, a direct proof can be found 
in \cite{Merkle2014}. Our proof is a consequence of Theorem 
\ref{thm:composition_comp_mon_1}.

Theorem \ref{thm:composition_comp_mon_2} is similar to Theorem 
\ref{thm:composition_comp_mon_1}. It also 
can be proved using the formula for the $n$-th derivative of a composite 
function, but the use of Leibniz's formula is enough. 

A direct proof of Corollary \ref{cor:log_f_diff_comp_mon} was given in 
\cite{QiChen2004}, here in Corollary \ref{cor:logcomp_comp}
 we obtain easily from Theorem 
\ref{thm:composition_comp_mon_2}.

Lemma \ref{lem:phi_x_per_x} was stated in \cite{Chen2007} as Lemma 1, but 
the proof of formula for $f^{(n)}(0)$ (\cite{Chen2007}, (5)) is wrong. The author 
says 

"It is known [3, p.26]
\begin{equation*}
    \sum_{k=0}^n\frac{(-1)^k\binom{n}{k}}{k+1}=\frac{1}{n+1}.
\end{equation*}
Using L'Hospital's rule, we have
\begin{align*}
    f^{(n)}(0) &=\sum_{k=0}^n \binom{n}{k} (-1)^k k!\lim_{x\to 0}
    \frac{\phi^{(n-k)}(x)}{x^{k+1}} \\
    &=\phi^{(n+1)}(0)\sum_{k=0}^n 
    \frac{(-1)^k\binom{n}{k}}{k+1}=\frac{1}{n+1}\phi^{(n+1)}(0).\text{"}
\end{align*}
The problem is that to apply L'Hospital's rule to calculate 
$\lim_{x\to 0}\frac{\phi^{(n-k)}(x)}
{x^{k+1}}$ we should assume that $\phi^{(j)}(0)=0$ ($j=0,\ldots,n$) 
because the denominator in the fraction is $0$. 

However using the correctly proved identity about $f^{(n)}(x)$ $(x\neq 0)$ 
(\cite{Chen2007} Lemma 1) 
 we can finish the proof easily
and we do not need to use the identity [3, p.26].  

A nice application of Lemma \ref{lem:phi_x_per_x} is Theorem 
\ref{thm:f_powerto_g_log_comp} where we prove the logarithmically  
complete monotonicity of a general function. This result will be used 
in the proof of  
Theorem \ref{thm:frac_gamma_power}.

Next result, Lemma \ref{lemma:1perx_comp_mon}, was proved in 
\cite{MillerSamko2001} Lemma 1 using the 
integral representation for $0<\mu<1$
\begin{equation*}
\left(1+\frac{1}{x}\right)^{\mu}=1+
\frac{\mu}{x^{\mu}}\int_1^{\infty}\frac{1}{t^{\mu+1}(xt+1)^{1-\mu}}\,dt.
\end{equation*}
We give a more general and simpler proof for $\mu>0$ not using the integral 
representation.

Theorem \ref{thm:linfraction_log_comp_mon} is about the 
 complete monotonicity of logarithm of linear fractional functions. 
Although it seems very elementary I was unable to find a reference in the 
literature. 

Theorem \ref{thm:fractionpower_log_comp_mon} is a generalization and an 
improvement of Lemma \ref{lemma:1perx_comp_mon}.

There are many results where the $\Gamma$ and the digamma, 
$\psi(z):=\Gamma'(z)/\Gamma(z)$ $(z\neq 0,-1,-2,\ldots)$, functions appear. 
Theorem \ref{thm:fraction_of_Gamma_log_comp} shows that the logarithm of 
the fraction of $\Gamma$ 
functions in general is not completely monotone. 

In \cite{Chen2007} Theorem 1 the following is stated.
Let $a\geq 0$, $b>0$ be given real numbers with $0<1-b+a<1$. Then, for all 
real numbers $\alpha$, the function $f_{\alpha}(x)=(x+a)^{\alpha}\left[ 
\psi(x+b)-\psi(x+a)-\frac{b-a}{x+a}\right]$ is strictly completely monotonic 
on $(-a,\infty)$ if and only if $\alpha\leq 1$. The proof of the "if" part 
is correct, but the proof 
of the "only if $\alpha\leq 1$" part is incorrect because it was 
calculated and used later 
\begin{align*}
	&f_{\alpha}'(x)= \\
    &(x+a)^{\alpha-1}\left\{
	\alpha\left[\psi(x+b)-\psi(x+a)\right]+
	(x+a)\left[\psi '(x+b)-\psi '(x+a)\right]
	\right\},
\end{align*}
which is wrong. 
We cannot prove the "only if" part, but we can give a necessary condition 
$\alpha\leq 2$ in Theorem \ref{thm:psi_Chen2007_improve_comp_mon} $(iii_2)$. 
Furthermore, if $-1-2a+2b>0$ then $\alpha<2$, Remark \ref{rem:alpha_less_2}. 
This motivates the following 
\begin{openproblem}
What is the value of $\alpha_0$ with the property that 
the function $f_{\alpha}(x)$ (see above) is strictly completely monotonic 
on $(-a,\infty)$ if and only if $\alpha\leq\alpha_0$?
\end{openproblem}
From \cite{Chen2007} Theorem 1 and Theorem \ref{thm:psi_Chen2007_improve_comp_mon} 
$(iii_2)$ we know that $1\leq\alpha_0\leq 2$. 

Theorem \ref{thm:psi_Chen2007_improve_comp_mon} $(iii_1)$ and $(iii_2)$ is a generalization 
of \cite{Chen2007} Theorem 1 wherein we can prove an 
"if and only if $\alpha\leq 1$" statement when $\beta<b-a$. 
Using our general Theorem 
\ref{thm:psi_Chen2007_improve_comp_mon} we will be able to prove 
the "if part" of 
Theorem \ref{thm:wrong_mon} in one row in a more general circumstance. 

In \cite{Chen2007} Theorem 2 it was proved that the function $g_{\beta}$ 
(see Theorem \ref{thm:wrong_mon}) is strictly logarithmically completely 
monotone on 
$(-a,\infty)$, reducing the problem to $g_{b-a}(x)$ by the identity 
$g_{\beta}(x)=(x+a)^{\beta-(b-a)}g_{b-a}(x)$.  We  
show that it can be proved directly without this identity.

The special case $b=1$, $a=1/2$, $\beta=0$, $c_0=\sqrt{\pi}$, 
$g(x)=1$ of Theorem \ref{thm:frac_gamma_power} 
was proved in \cite{Chen2007} Theorem 3.

Theorem \ref{thm:log_gamma_vogt_general} in the special case 
$\beta=0$ was proved in \cite{QiChen2004} Theorem 2. 
Our more general result is an easy consequence of 
Lemma \ref{lem:vogt_comp_mon} that was proved in \cite{VogtVoigt2002} Theorem 1, 
and Theorem \ref{thm:linfraction_log_comp_mon}.

\section{Lemmas and Theorems}

\begin{theorem}\label{thm:analytic_abs}
Let $I:=(a,b)$ be finite interval, $f$ be absolutely monotone in $I$. Then 
$f$ can be extended analytically into the complex plane ($z=x+iy$), and the 
extended function $f(z)$ will be analytic in the circle $|z-a|<b-a$, 
and will be absolutely monotone in $[a,b)$. 

If $I=(-\infty,b)$ and $b$ is finite, then $f(z)$ will be analytic in the 
half-plane $\mathrm{Re}\,z<b$. If $a$ is finite and $I=(a,\infty)$, then $f(z)$ 
will be analytic in the whole plane. If $I=(-\infty,\infty)$ 
then $f(z)$ will be analytic in the whole plane.
\end{theorem}

\begin{corollary}\label{cor:analytic_comp}
Let $I:=(a,b)$ be finite interval, $g$ be completely monotone on $I$. 
Then it can be extended analytically into the complex $z$-plane 
($z=x+iy$), and the extended function $g(z)$ will be analytic in the circle 
$|z-b|<b-a$, and will be completely monotone in $(a,b]$.

If $I=(-\infty,b)$ and $b$ is finite, then $g(z)$ will be analytic in the 
whole plane. If $a$ is finite and $I=(a,\infty)$, then 
$g(z)$ 
will be analytic in the half-plane $\mathrm{Re}\,z>a$. If $I=(-\infty,\infty)$ 
then $g(z)$ will be analytic in the whole plane.
\end{corollary}

\begin{theorem}\label{thm:product_of_comp_mon}
Let $I\subseteq\mathbf{R}$ be interval, $f,g:I\to\mathbf{R}$ be arbitrary 
completely monotone functions. Then $fg$ is completely monotone.
\end{theorem}

\begin{corollary}\label{cor:f_power_g_logcomp_mon}
Let $I\subseteq\mathbf{R}$ be an interval, $f:I\to\mathbf{R}$,  
$\log\circ f$ be a
 completely monotone function, $g:I\to\mathbf{R}$ be a
completely monotone function. 
Then $f^g$ is logarithmically completely monotone.
\end{corollary}

\begin{theorem}\label{thm:composition_comp_mon_1}
Let $I\subseteq\mathbf{R}$ be an interval, $g:I\to\mathbf{R}$, $g\geq 0$, $g'$ be 
completely monotone, $f$ be completely monotone in the set
$\mathrm{Range}\,g$. Then $f\circ g$ is completely monotone.
\end{theorem}

\begin{corollary}\label{cor:f_x_alpha}
If $0<\alpha\leq 1$,  
$f$ is completely monotone in $(0,b)$, $b$ is finite or infinite,
then $f(x^{\alpha})$ 
is completely monotone in $(0,b^{1/\alpha})$.
\end{corollary}

\begin{corollary}\label{cor:log_f_diff_comp_mon}
Let $f:(a,b)\to\mathbf{R}$, $f>0$. If $f$ is logarithmically completely 
monotone in $(a,b)$ then $f$ is completely monotone in $(a,b)$.
\end{corollary}

\begin{theorem}\label{thm:composition_comp_mon_2}
Let $I\subseteq\mathbf{R}$ be interval, $g:I\to\mathbf{R}$, $-g'$    
be completely monotone, $f$ be absolutely monotone in the set 
$\mathrm{Range}\,g$. Then $f\circ g$ is completely monotone.
\end{theorem}

\begin{corollary}\label{cor:logcomp_comp}
Let $h:(a,b)\to\mathbf{R}$, $h>0$. If $h$ is logarithmically 
completely monotone 
in $(a,b)$ then $h$ is completely monotone in $(a,b)$.
\end{corollary}

\begin{lemma}\label{lem:phi_x_per_x}
Let the function $\varphi$ have derivatives of all orders on $[a,b)$, 
$a\leq 0<b$, and 
$\varphi(0)=0$. Define the function $f$ by
\begin{equation*}
	f(x)=
	\begin{cases}
		\frac{\varphi(x)}{x}, & x\in [a,b)-\{0\},\\
		\varphi'(0), & x=0.
	\end{cases}
\end{equation*}
Then
\begin{equation*}
	f^{(n)}(x)=
	\begin{cases}
		\frac{1}{x^{n+1}}\sum_{k=0}^n \binom{n}{k}(-1)^k k!x^{n-k}\varphi^{(n-k)}(x),
		& x\neq 0,\\
		\frac{1}{n+1}\varphi^{(n+1)}(0), & x=0,
	\end{cases}
\end{equation*}
and $f^{(n)}(x)$ is continuous on $[a,b)$. 
Moreover 
\begin{equation*}
	\frac{\mathrm{d}}{\mathrm{d}x}\sum_{k=0}^n \binom{n}{k}(-1)^k k!x^{n-k}
	\varphi^{(n-k)}(x)=x^n\varphi^{(n+1)}(x).
\end{equation*}
\end{lemma}

\begin{theorem}\label{thm:f_powerto_g_log_comp}
Let  $f:[0,b)\to\mathbf{R}$,  assume  
$\log\circ f$ is a Bernstein function on $[0,b)$, 
$g$ is completely monotone on $(0,b)$, then $f(x)^{g(x)/x}$ 
is logarithmically completely monotone on $(0,b)$.
\end{theorem}

\begin{lemma}\label{lemma:1perx_comp_mon}
The function 
\begin{equation*}
f(x):=\left(a+\frac{b}{x}\right)^{\mu}\quad(a,b,\mu\geq 0)
\end{equation*}
is completely monotone in $(0,\infty)$. 
\end{lemma}

\begin{theorem}\label{thm:linfraction_log_comp_mon}
Let $I:=(\alpha,\infty)$ be interval, 
$a,c\geq 0$, 
$a^2+b^2, c^2+d^2\neq 0$. Then 
\begin{equation*}
	f(x):=\log\left(\frac{ax+b}{cx+d}\right)
\end{equation*}
is completely monotone in $I$ if and only if one of the following conditions holds 
\begingroup
\renewcommand\labelenumi{(\theenumi)}
\begin{enumerate} 
	\item $a=0$, $c=0$, $0<d\leq b$ or $0>d\geq b$, $\alpha$ is arbitrary;
	\item $a>0$, $c=a$, $d=b$, $\alpha\geq -d/c$;
	\item $a>0$, $c=a$, $b>d$, $\alpha\geq -d/c$;
	\item $a>0$, $a>c>0$, $ad-bc\leq 0$, $\alpha\geq -d/c$.
\end{enumerate}
\endgroup
\end{theorem}

\begin{remark}
	If $s+A,s+B>0$ then
	\begin{equation*}
		\log\left(\frac{s+A}{s+B}\right)=\int_0^{\infty}e^{-st}\,
		\frac{e^{-Bt}-e^{-At}}{t}\,dt.
	\end{equation*}
\end{remark}

\begin{theorem}\label{thm:fractionpower_log_comp_mon}
Let $g$ be completely monotone function in $(x_0,\infty)$, and 
\\
$\alpha:=\max(x_0,-d/c)$. If

\textup{(1)} $a>0$, $c=a$, $b\geq d$;

or

\textup{(2)} $a>0$, $a>c>0$, $ad-bc\leq 0$,

then 
\begin{equation*}
	f(x):=\left(\frac{ax+b}{cx+d}\right)^{g(x)}
\end{equation*}
is logarithmically completely monotone, and completely monotone 
in $(\alpha,\infty)$.
\end{theorem}

\begin{theorem}\label{thm:fraction_of_Gamma_log_comp}
Let $a,c>0$, and 
\begin{equation*}
	f(x):=\log\left(\frac{\Gamma(ax+b)}{\Gamma(cx+d)}\right),
\end{equation*}
where $x\in (\alpha,\infty)=:I$, $\alpha:=\max\{-b/a,-d/c\}$. Then $f$ is 
completely monotone in $I$ if and only if $a=c$, $b=d$.
\end{theorem}

\begin{theorem}\label{thm:psi_Chen2007_improve_comp_mon}
Let $a,b\geq 0$, $\alpha,\beta$ be given real numbers and define the function 
$f_{\alpha,\beta}(x)=(x+a)^{\alpha} [\psi(x+b)-\psi(x+a)-\frac{\beta}{x+a}]$ 
for $x>\max (-a,-b)$. 

$(i)$ If $b=a$ then $f_{\alpha,\beta}(x)$ is completely monotone on 
$(-a,\infty)$ if and only if one of the following conditions holds

$(i_1)$ $\beta=0$, $\alpha$ is arbitrary;

$(i_2)$ $\beta<0$, $\alpha\leq 1$.

(ii) If $b=a+1$ then $f_{\alpha,\beta}(x)$ is completely monotone on 
$(-a,\infty)$ if and only if one of the following conditions holds 

$(ii_1)$ $\beta=1$, $\alpha$ is arbitrary;

$(ii_2)$ $\beta<1$, $\alpha\leq 1$.

(iii) Let $b\neq a, a+1$. 
 Assume that $0<1-b+a<1$, and 
$\beta\leq b-a$. 

$(iii_1)$ Let $\alpha\leq 1$. 
Then the function $f_{\alpha,\beta}(x)$ is completely 
monotone on $(-a,\infty)$.  
The lower and upper bounds for $b$, $a<b$ and $b<a+1$, are sharp, and the 
upper bound for $\beta$ is sharp.

$(iii_2)$ Let the function $f_{\alpha,\beta}(x)$ be completely 
monotone on $(-a,\infty)$. If $\beta<b-a$ then $\alpha\leq 1$, if $\beta=b-a$ 
then $\alpha\leq 2$.
\end{theorem}

\begin{remark}\label{rem:alpha_less_2}
If $\beta=b-a$ and $-1-2a+2b>0$ in Theorem 
\ref{thm:psi_Chen2007_improve_comp_mon} $(iii_2)$ 
then $\alpha<2$. 
\end{remark}

\begin{theorem}\label{thm:wrong_mon}
Let $a\geq 0$, $b>0$ be given real numbers with $0<1-b+a<1$. 
 Then the function 
$g_{\beta}(x)=\frac{(x+a)^{\beta}\Gamma(x+a)}{\Gamma(x+b)}$ is  
logarithmically completely monotone on $(-a,\infty)$  
 if and only if $\beta\leq b-a$.
\end{theorem}

\begin{theorem}\label{thm:frac_gamma_power}
Let $a,b>0$, $\beta$ be given real numbers with 
$0<1-b+a\leq 1$, $\beta\leq b-a$, $c_0\geq 
a^{\beta}\Gamma(a)/\Gamma(b)$, $g(x)$ be completely monotone on 
$(0,\infty)$. Then the function 
$\left[\frac{c_0\Gamma(x+b)}{(x+a)^{\beta}\Gamma(x+a)}\right]^{g(x)/x}$ 
is logarithmically completely monotone 
on $(0,\infty)$.
\end{theorem}

\begin{lemma}\label{lem:vogt_comp_mon}
The function $1+\frac{1}{x}\log\,\Gamma(x+1)-\log(x+1)$ 
is completely monotone on $(-1,\infty)$.
\end{lemma}

\begin{theorem}\label{thm:log_gamma_vogt_general}
Let $0\leq\beta \leq 1$. Then the function $1+\frac{1}{x}\log\,\Gamma(x+1)-\log(x+\beta)$ 
is completely monotone on $(0,\infty)$.
\end{theorem}

\section{Proofs}

\textbf{Proof of Theorem \ref{thm:analytic_abs}.} 
We follow \cite{Widder1946} with a little modification. 
Let $u:=a+\varepsilon$, where $0<\varepsilon<(b-a)/3$. Then $f$ is absolutely monotone in $u\leq x<b$. 
By Taylor's formula for every $u\leq x<b$ we have
\begin{equation*}
	f(x)=\sum_{k=0}^{n} \frac{f^{(k)}(u)}{k!} (x-u)^k+R_n(x), 
\end{equation*}
where
\begin{align*}
	R_n(x):&=\int_{u}^x \frac{f^{(n+1)}(t)}{n!}(x-t)^{n}\,dt\\
	&=\frac{(x-u)^{n+1}}{n!}\int_0^1 f^{(n+1)}(u+[x-u]t)(1-t)^n\,dt.
\end{align*}
Since $f^{(n+2)}\geq 0$ therefore $f^{(n+1)}(u+[x-u]t)$ is a monotone 
increasing function of $x$ when $t$ is fixed, so that if $u\leq x\leq c<b$
we have
\begin{align*}
	0 &\leq R_n(x)\leq\frac{(x-u)^{n+1}}{n!}\int_0^1 (1-t)^n  
	f^{(n+1)}(u+[c-u]t)\,dt\\
	&=\frac{(x-u)^{n+1}}{(c-u)^{n+1}}
	\left[
	f(c)-f(u)-f'(u)(c-u)-\ldots -f^{(n)}(u)\frac{(c-u)^n}{n!}
	\right] \\
	&\leq f(c)\left(\frac{x-u}{c-u}\right)^{n+1}.
\end{align*}
Hence 
\begin{equation*}
	\lim_{n\to\infty}R_n(x)=0\quad(u\leq x<c<b).
\end{equation*}
Since $c$ is arbitrary we obtain
\begin{equation*}
	f(x)=\sum_{n=0}^{\infty}\frac{f^{(n)}(u)}{n!}(x-u)^n \quad(u\leq x<b).
\end{equation*}
The analytical extension is 
\begin{equation*}
	\widetilde{f}(z):=\sum_{n=0}^{\infty}\frac{f^{(n)}(u)}{n!}(z-u)^n,
\end{equation*}
and $\widetilde{f}$ will be analytic in the circle $|z-u|<b-u$. 
Since $\varepsilon>0$ was arbitrary the statement follows. $\qedsymbol$

\textbf{Proof of Corollary \ref{cor:analytic_comp}.} 
Define the function $f(x):=g(-x)$. Then $f$ is absolutely monotone 
in $(-b,-a)$ so we can apply Theorem \ref{thm:analytic_abs}. $\qedsymbol$

\textbf{Proof of Theorem \ref{thm:product_of_comp_mon}.} 
We prove by induction. 
For $n=0,1$ it is easy. Our induction hypothesis is that if $u,v$ are any 
completely monotone functions then $(-1)^k(uv)^{(k)}\geq 0$ for $k\leq n$. 
Now 
\begin{align*}
	(-1)^{n+1}(fg)^{(n+1)} &=(-1)^n((-f')g+f(-g'))^{(n)} \\
	&=	(-1)^n((-f')g)^{(n)}+(-1)^n(f(-g'))^{(n)}.
\end{align*}
Since $-f', -g'$ are completely monotone the induction hypothesis gives the 
non-negativity. $\qedsymbol$

\textbf{Proof of Corollary \ref{cor:f_power_g_logcomp_mon}.}
We know
\begin{equation*}
	\log\left(f^g\right)=g\cdot \log\circ f.
\end{equation*}
Since $\log\circ f$ is completely monotone,  
Theorem \ref{thm:product_of_comp_mon} yields the statement. $\qedsymbol$

\textbf{Proof of Theorem \ref{thm:composition_comp_mon_1}.} 
We prove by induction. If $n=0$ then 

$(-1)^n(f\circ g)^{(n)}(x)=f(g(x))\geq 0$. 
If $n=1$ then 
\begin{equation*}
	(-1)(f\circ g)'=(-f')\circ g\cdot g'\geq 0.
\end{equation*}
Using Leibniz's formula 
we have
\begin{align*}
	(-1)^{n+1}(f\circ g)^{(n+1)} &=(-1)^{(n+1)}[f'\circ g \cdot g']^{(n)} \\
	&=(-1)^n \sum_{k=0}^n \binom{n}{k}((-f')\circ g)^{(k)}\cdot (g')^{(n-k)} \\
	&=\sum_{k=0}^n \binom{n}{k}(-1)^k((-f')\circ g)^{(k)}\cdot 
	(-1)^{(n-k)}(g')^{(n-k)}.
\end{align*}
Here $-f'$ is completely monotone. 
Assuming the validity of induction hypothesis for $n$, for any completely 
monotone function, and for any non-negative function whose first derivative is 
completely monotone, we obtain the statement. $\qedsymbol$

\textbf{Proof of Corollary \ref{cor:f_x_alpha}.}
If $\alpha=1$ the statement is obvious, so we assume that $0<\alpha<1$. 
Denote $g(x):=x^{\alpha}$. Then $g\geq 0$ and $g'$ is completely monotone, 
thus Theorem \ref{thm:composition_comp_mon_1} yields the statement.
 $\qedsymbol$

\textbf{Proof of Corollary \ref{cor:log_f_diff_comp_mon}.}
We choose $\delta_1,\delta_2$ such that $a<\delta_1<\delta_2<b$.  
Denote $M:=\max\,\{f(x)\,|\,\delta_1\leq x\leq\delta_2\}$. 
Denote $g(x):=-(\log\circ (f/(2M)))(x)$, $\delta_1\leq x\leq \delta_2$. 
Then $g\geq 0$, and 
$g'(x)=(-\log\circ f)'$ is completely monotone in $(\delta_1,\delta_2)$. 
Take $h(x):=e^{-x}$. Then $h$ is completely monotone in $(0,\infty)$. 
By Theorem \ref{thm:composition_comp_mon_1} 
$(h\circ g)(x)=f(x)/(2M)$ 
is completely monotone in $(\delta_1,\delta_2)$. Since 
$(\delta_1,\delta_2)$ 
is arbitrary subinterval of $(a,b)$, and $f(x)=(2M)(f(x)/(2M))$, the 
corollary follows. $\qedsymbol$

\textbf{Proof of Theorem \ref{thm:composition_comp_mon_2}.}
We prove by induction.  If $n=0$ then 

$(-1)^n(f\circ g)^{(n)}(x)=f(g(x))\geq 0$. 
If $n=1$ then 
\begin{equation*}
	(-1)(f\circ g)'=f'\circ g\cdot (-g')\geq 0.
\end{equation*}
Using Leibniz's formula we have 
\begin{align*}
	(-1)^{n+1}(f\circ g)^{(n+1)} &=(-1)^{n}[f'\circ g \cdot
	(-g')]^{(n)} \\
	&=(-1)^n \sum_{k=0}^n \binom{n}{k}(f'\circ g)^{(k)}\cdot 
	(-g')^{(n-k)} \\
	&=\sum_{k=0}^n \binom{n}{k}(-1)^k(f'\circ g)^{(k)}\cdot 
	(-1)^{(n-k)}(-g')^{(n-k)}.
\end{align*}
Here $-g'$ is completely monotone, $f'$ is absolutely monotone. 
Assuming the validity of induction hypothesis for $n$, for any 
absolutely monotone function, and for any function $g$ such that $-g'$ is 
completely monotone,
we obtain the statement. $\qedsymbol$

\textbf{Proof of Corollary \ref{cor:logcomp_comp}.}
Take $f:=\exp$, and $g:=\log\circ h$. $\qedsymbol$

\textbf{Proof of Lemma \ref{lem:phi_x_per_x}.} We prove only $f^{(n)}(0)=
\frac{1}{n+1}\varphi^{(n+1)}(0)$ and the continuity of $f^{(n)}$, because 
everything else is proved correctly in 
\cite{Chen2007} Lemma 1. 
For $n=0$ this and the continuity follow from 
the definition of $f$, 
so we may assume that $n\geq 1$. We prove by induction. 
\begin{align*}
    f^{(n)}(0) &=\lim_{x\to 0} \frac{f^{(n-1)}(x)-f^{(n-1)}(0)}{x-0}\\
    &=\lim_{x\to 0} f^{(n)}(x),
\end{align*}
assuming that the assumptions of L'Hospital's rule are satisfied. 
Since $f^{(n-1)}$ is continuous the fraction has the type "$0/0$" 
at $x=0$. 
We know from \cite{Chen2007} Lemma 1 that 
\begin{equation*}
    f^{(n)}(x)=\frac{1}{x^{n+1}}\sum_{k=0}^n\binom{n}{k}(-1)^k k! x^{n-k}\varphi^{(n-k)}(x),\quad x\neq 0.
\end{equation*}
The sum is $0$ at $x=0$ because of $\varphi(0)=0$, so we use 
L'Hospital's rule 
to determine $\lim_{x\to 0}f^{(n)}(x)$.
\begin{align*}
    \lim_{x\to 0} f^{(n)}(x) &=\lim_{x\to 0}
    \frac{\sum_{k=0}^n\binom{n}{k}(-1)^k k! x^{n-k}\varphi^{(n-k)}(x)}
    {x^{n+1}} \\
    &=\lim_{x\to 0}
    \frac{\frac{\mathrm{d}}{\mathrm{d}x}\sum_{k=0}^n \binom{n}{k}(-1)^k k!x^{n-k}
	\varphi^{(n-k)}(x)} {(n+1)x^n}.
\end{align*}
In the last fraction the numerator is $x^n\varphi^{(n+1)}(x)$ that was proved for 
any $n$ and $x\in [a,b)$. Since $\varphi^{(n+1)}(x)$ is continuous on $[a,b)$ we have 
$f^{(n)}(0)=\frac{1}{n+1}\varphi^{(n+1)}(0)$. $\qedsymbol$

\textbf{Proof of Theorem \ref{thm:f_powerto_g_log_comp}.}
By Corollary \ref{cor:f_power_g_logcomp_mon} it is enough to 
prove the case $g=1$. We follow the idea of 
\cite{Chen2007} Proof of Theorem 3. 
Define the function
\begin{equation*}
	\varphi(x):=(\log\circ f)(x). 
\end{equation*}
Since $\varphi$ is a Bernstein function thus   $\varphi(0)\geq 0$. 
By Lemma \ref{lem:phi_x_per_x} we obtain on the 
interval $(0,b)$ for $n\geq 1$ 
\begin{align*}
	&x^{n+1}\left(\frac{1}{x}\varphi(x)\right)^{(n)}=
	\sum_{k=0}^n \binom{n}{k}(-1)^k k!x^{n-k}\varphi^{(n-k)}(x)=:\Phi(x) \\
	&\Phi '(x)=x^n \varphi^{(n+1)}(x).
\end{align*}
Since the function $\varphi$ has derivatives of all orders on $[0,b)$  
and $\varphi(0)\geq 0$, it easily follows that $(-1)^n\Phi(0)\geq 0$. 

If $n$ is even, then we have for $x>0$ that 
$(\varphi')^{(n)}(x)\geq 0$ which implies $\Phi'(x)\geq 0$, so 
$\Phi(x)\geq 0$, and consequently $(\varphi(x)/x)^{(n)}=
(-1)^n (\varphi(x)/x)^{(n)} \geq 0$. 

Similarly, if $n$ is odd, then we have for $x>0$ that 
$(\varphi')^{(n)}(x)\leq 0$ which implies $\Phi'(x)\leq 0$, so 
$\Phi(x)\leq 0$, and consequently $(\varphi(x)/x)^{(n)}\leq 0$, and 
it gives $(-1)^n (\varphi(x)/x)^{(n)} \geq 0$. 

Hence 
\begin{equation*}
	(-1)^n (\varphi(x)/x)^{(n)} \geq 0\quad(x> 0).
\end{equation*}
The theorem is proved.  $\qedsymbol$

\textbf{Proof of Lemma \ref{lemma:1perx_comp_mon}.}
The cases $b=0$, or $\mu=0$ are obvious. If $a=0$  then 
\begin{equation*}
	(-1)^n\left(\frac{1}{x^{\mu}}\right)^{(n)}=\mu(\mu+1)\cdot 
	(\mu+n-1) x^{-(\mu+n)}>0.
\end{equation*}
So assume that $a,b,\mu>0$. Then
\begin{equation*}
	f(x):=\left(a+\frac{b}{x}\right)^{\mu}=
	a^{\mu}\left(1+\frac{b/a}{x}\right)^{\mu}=
	a^{\mu}\left(1+\frac{1}{ax/b}\right)^{\mu}. 
\end{equation*}
It gives
\begin{equation*}
	(-\log(f(x)))'=\mu\left(\frac{1}{x}-\frac{1}{x+1/\alpha}\right),
\end{equation*}
where $\alpha=a/b$. Since 
\begin{equation*}
	(-1)^n\left(\frac{1}{x+c}\right)^{(n)}=n!(x+c)^{-(n+1)}
\end{equation*}
we obtain
\begin{equation*}
	(-1)^n\mu\left(\frac{1}{x}-\frac{1}{x+1/\alpha}\right)^{(n)}=
	\mu n!\left(\frac{1}{x^{n+1}} -\frac{1}{(x+1/\alpha)^{n+1}} \right)>0.
\end{equation*}
Applying Corollary \ref{cor:log_f_diff_comp_mon} we get that $f(x)$ is 
completely monotone in $(0,\infty)$. $\qedsymbol$

\textbf{Proof of Theorem \ref{thm:linfraction_log_comp_mon}.}
\begingroup
\renewcommand\labelenumi{(\theenumi)}
First we consider the "only if" part, then the "if" part is an easy verification. 
\begin{enumerate} 
	\item $a=0$, $b>0$. The $c>0$ case is impossible because $f\geq 0$, so 
	$c=0$. Then $0<d\leq b$, and $f$ is a constant function.
	\item $a>0$. Then $c>0$. Since $\lim_{x\to\infty}f(x)=\log(a/c)$ and 
	$f\geq 0$ we have $a\geq c$. 
	If $a=c$, and $d=b$, then $f$ is a constant function. 
	\item If $d\neq b$, then since $f\geq 0$ we get $b>d$. Now we can write
	\begin{align*}
		f(x) &=\log\left(\frac{ax+b}{ax+d}\right) \\[1.1em]
		&=\log(x+b/a)-\log(x+d/a).
	\end{align*}
	From this it follows
	\begin{equation*}
		f'(x)=\frac{1}{x+b/a}-\frac{1}{x+d/a}<0.
	\end{equation*}
	If $n\geq 2$ then
	\begin{equation*}
		(-1)^nf^{(n)}(x)=-(n-1)!\left(\frac{1}{(x+b/a)^n}-
		\frac{1}{(x+d/a)^n}\right)>0.
	\end{equation*}
	\item $a>c$. Then $c>0$. The condition $f\geq 0$ implies $ax+b\geq cx+d$ 
	from which we obtain
	\begin{equation}\label{eq:alpha_value}
		x\geq \frac{d-b}{a-c}.
	\end{equation}
	We can write
	\begin{equation*}
		f(x)=\log(a/c)+\log(x+b/a)-\log(x+d/c).
	\end{equation*}
	This yields
	\begin{equation*}
		f'(x)=\frac{1}{x+b/a}-\frac{1}{x+d/c}=\frac{d/c-b/a}{(x+b/a)(x+d/c)}.
	\end{equation*}
	The condition $f'\leq 0$ gives $ad-bc\leq 0$. This inequality shows that 
	in \eqref{eq:alpha_value} $\frac{d-b}{a-c}\leq -d/c$, hence $\alpha\geq -d/c$. 
	If $n\geq 2$ then
	\begin{equation*}
		(-1)^nf^{(n)}(x)=-(n-1)!\left(\frac{1}{(x+b/a)^n}-\frac{1}{(x+d/c)^n}\right)
		\geq 0.
	\end{equation*}
	By Corollary \ref{cor:logcomp_comp} the proposition 
	follows.\quad$\qedsymbol$
\end{enumerate}
\endgroup

\textbf{Proof of Theorem \ref{thm:fractionpower_log_comp_mon}.}
Obviously
\begin{equation*}
	\log(f(x))=g(x)\log\left(\frac{ax+b}{cx+d}\right).
\end{equation*}
By Theorem \ref{thm:linfraction_log_comp_mon}, Theorem 
\ref{thm:product_of_comp_mon}, and Corollary \ref{cor:logcomp_comp} 
the statement follows. $\qedsymbol$

\textbf{Proof of Theorem \ref{thm:fraction_of_Gamma_log_comp}.}
We know that \cite{Bateman1953} 1.18 (1)
\begin{equation*}
	\log\,\Gamma(z)=(z-1/2)\log\,z-z+1/2\,\log(2\pi)+O(1/z),\quad
	|\arg \,z|<\pi-\varepsilon,\,\varepsilon>0.
\end{equation*}
Since $f\geq 0$ we need to have $\log\,\Gamma(ax+b)\geq\log\Gamma(cx+d)$ for 
large $x>0$ also. So we get the condition
\begin{equation*}
	ax\log\,(ax+b)+O(x)\geq cx\log\,(cx+d)+O(x),
\end{equation*}
from which 
\begin{equation*}
	a/c\geq \log\,(cx+d)/\log\,(ax+b)+O(1/\log\,x).
\end{equation*}
Since the limit of the right-hand-side is $1$ we obtain $a\geq c$. 

The function $f$ is monotone decreasing thus for $x\geq x_0>\alpha$ it is bounded 
above with, say $c_0>0$.
Using again the asymptotics we have $a\leq c$. Hence $a=c$. Since we want 
$f\geq 0$ in the interval, and $\Gamma$ is strictly monotone increasing 
on $[1,\infty)$ we need to assume that $b\geq d$. 

We know that \cite{Bateman1953} 1.18 (4)
\begin{equation*}
	\frac{\Gamma(z+\alpha)}{\Gamma(z+\beta)}=z^{\alpha-\beta}(1+O(1/z)) \quad
	|\arg\,z|<\pi-\varepsilon,\,\varepsilon>0.
\end{equation*} 
Applying this we have
\begin{equation*}
	f(x)=\log\,a^{b-d}+(b-d)\log\,x+O(1/x).
\end{equation*}
The function $f$ is monotone decreasing thus for $x\geq x_0>\alpha$ it is bounded 
above with, say $c_0>0$. It implies $b\leq d$.  Hence $b=d$. 
 $\qedsymbol$

\textbf{Proof of Theorem \ref{thm:psi_Chen2007_improve_comp_mon}.}

(i) In this case $f_{\alpha,\beta}(x)=-\beta (x+a)^{\alpha-1}$. If $\beta=0$ 
then $f_{\alpha,\beta}(x)$ is a non-negative constant function, so completely 
monotone. If $\beta\neq 0$ then it is necessary that $\beta<0=b-a$. In this case 
$\alpha\leq 1$ is the necessary and sufficient condition for $f_{\alpha,\beta}(x)$ 
to be completely monotone. 

(ii) Since 
\begin{equation*}
    \psi(z+1)=\psi(z)+\frac{1}{z},
\end{equation*}
we obtain
\begin{equation*}
    f_{\alpha,\beta}(x)=(1-\beta)(x+a)^{\alpha-1}.
\end{equation*}
If $\beta=1$ then $f_{\alpha,\beta}(x)$ is a non-negative constant function, so 
completely monotone. If $\beta\neq 1$ then it is necessary that $\beta<1=b-a$. 
In this case 
$\alpha\leq 1$ is the necessary and sufficient condition for $f_{\alpha,\beta}(x)$ 
to be completely monotone. 

$(iii_1)$
To prove that $f_{\alpha,\beta}$ is completely monotone, we need to modify 
the function $\varphi$ in \cite{Chen2007} Proof of Theorem 1, 
\begin{equation*}
	\varphi(t):=1-\beta-\frac{e^{(1-b+a)t}-1}{e^t-1},
\end{equation*}
and the proof works. 

We show that the 
lower and upper bounds for $b$, $a<b$ and $b<a+1$, are sharp. 
Assume that $f_{1}:=f_{1,b-a}$ is completely monotone. 

We know that \cite{Bateman1953} 1.7 (3)
\begin{equation*}
	\psi(z)=-\gamma-\frac{1}{z}+\sum_{k=1}^{\infty}\frac{z}{k(z+k)},
\end{equation*}
uniformly on compact sets not containing the poles of $\psi$ 
($z=0,-1,-2,\ldots$), where $\gamma$ is the Euler–Mascheroni constant. 
This equation can be written into
\begin{equation*}
	\psi(z)=-\gamma-\frac{1}{z}+\sum_{k=1}^{\infty}\left[\frac{1}{k}-
	\frac{1}{z+k}\right].
\end{equation*}
From this we obtain
\begin{align}
	f_{1}(x) &=(x+a)[\psi(x+b)-\psi(x+a)]-(b-a) \nonumber\\
	&=(x+a)\left(\frac{1}{x+a}-\frac{1}{x+b}+\sum_{k=1}^{\infty}
	\left[\frac{1}{x+a+k}-\frac{1}{x+b+k}\right]
	\right) -(b-a) \nonumber\\
	&=\frac{b-a}{x+b}+(x+a)\sum_{k=1}^{\infty}\frac{b-a}{(x+a+k)(x+b+k)}-(b-a)
	\nonumber\\ 
	&=(b-a)(F(a,b,x)-1), \label{eq:f1}
\end{align}
where
\begin{equation*}
	F(a,b,x)=\frac{1}{x+b}+(x+a)\sum_{k=1}^{\infty}\frac{1}{(x+a+k)(x+b+k)}.
\end{equation*}
The function $F$ is strictly monotone decreasing in $b$ (if $a,x$ is fixed),  
and  
\begin{align*}
	F(a,a+1,x) &=\frac{1}{x+a+1}+(x+a)\sum_{k=1}^{\infty}\left[
	\frac{1}{x+a+k}-\frac{1}{x+a+1+k}\right]\\
	&=\frac{1}{x+a+1}+(x+a)\frac{1}{x+a+1}=1.
\end{align*}
If $b>a$ then \eqref{eq:f1} implies $b<a+1$. If $b<a$ then we should have  
$b\geq a+1$, but this is impossible. 

Now we prove that the assumption $\beta\leq b-a$ is sharp. 
We know \cite{Luke1964} 2.11 (9)  that
\begin{equation*}
	\psi(z+a)=\log\,z+ \frac{a-1/2}{z}
	+O\left(\frac{1}{z^2}\right),
	\quad (z\to\infty\,\mathrm{in}\, |\arg\,z|<\pi-\varepsilon,\,\varepsilon>0),
\end{equation*}
where $a$ is an arbitrary complex number. 
From this
\begin{equation*}
	\psi(x+b)=\log\,x+\frac{b-1/2}{x}+
	O\left(\frac{1}{x^2}\right).
\end{equation*}
Thus we get
\begin{equation}\label{eq:psi_different}
	\psi(x+b)-\psi(x+a) =\frac{b-a}{x}
	+O\left(\frac{1}{x^2}\right).
\end{equation}
Hence we obtain
\begin{align*}
	f_{1,\beta}(x) &=(x+a)\left(\frac{b-a-\beta}{x}+
	O\left(\frac{1}{x^2}\right) \right)\\
	&=b-a-\beta+O\left(\frac{1}{x}\right).
\end{align*}
Since $f_{1,\beta}(x)\geq 0$ for all $x>-a$, we should assume that 
$\beta\leq b-a$.

$(iii_2)$
We follow the idea from \cite{Chen2007} Theorem 1. 
Assume that $f_{\alpha,\beta}$ is completely monotone on $(-a,\infty)$. 
Then we have for all $x>-a$
\begin{align*}
	f_{\alpha,\beta}'(x)&=
    (x+a)^{\alpha-2}\times\\
    &\left\{
	\alpha(x+a)[\psi(x+b)-\psi(x+a)] +(x+a)^2
	[\psi'(x+b)-\psi'(x+a)]\right. \\
    &\left. -(\alpha-1)\beta\right\} < 0,
\end{align*}
which implies
\begin{align}\label{eq:upper_for_alpha_beta}
	&\alpha\{ (x+a)[\psi(x+b)-\psi(x+a)]-\beta\} \nonumber\\
    &<-(x+a)^2[\psi'(x+b)-\psi'(x+a)]-\beta.
\end{align}
We know \cite{Luke1964} 2.11 (9)  that
\begin{align*}
	\psi(z+a) & =\log\,z- \sum_{k=0}^1
\frac{(-1)^{k+1}B_{k+1}(a)}{k+1}z^{-k-1}+O\left(\frac{1}{z^3}\right), \\
	& (z\to\infty\,\mathrm{in}\, |\arg\,z|<\pi-\varepsilon,\,\varepsilon>0),
\end{align*}
where $a$ is an arbitrary complex number, and $B_{k+1}(a)$ are Bernoulli polynomials, 
$B_1(a)=a-1/2$, $B_2(a)=a^2-a+1/6$, that is, 
\begin{equation*}
    \psi(x+a)=\log x+\frac{a-1/2}{x}-\frac{a^2-a+1/6}{2}\frac{1}{x^2}+
    O\left(\frac{1}{x^3}\right). 
\end{equation*}
Thus we get 
\begin{align}
    & (x+a)[\psi(x+b)-\psi(x+a)]-(b-a) \nonumber \\
    & =(x+a) \left[
    \frac{b-a}{x}+\frac{(a-b)(a+b-1)}{2}\frac{1}{x^2}+
    O\left(\frac{1}{x^3}\right) \right]-(b-a) \nonumber \\
    & =\frac{(b-a)(a-b+1)}{2}\frac{1}{x} +O\left(\frac{1}{x^2}\right).
    \label{eq:psi_difference_asymp_finer}
\end{align}
Since $\psi$ is analytical as $z\to\infty$ $(|\arg\,z|<\pi-\varepsilon,\,
\varepsilon>0)$,  
we obtain
\begin{align*}
	\psi'(z+a) & =\frac{1}{z}+ \sum_{k=0}^1
(-1)^{k+1}B_{k+1}(a)z^{-k-2}+O\left(\frac{1}{z^4}\right), \\
	& (z\to\infty\,\mathrm{in}\, |\arg\,z|<\pi-\varepsilon,\,\varepsilon>0),
\end{align*}
that is, 
\begin{equation*}
    \psi'(x+a)=\frac{1}{x}-\frac{a-1/2}{x^2}+\frac{a^2-a+1/6}{x^3}+
    O\left(\frac{1}{x^4}\right). 
\end{equation*}
Thus we get
\begin{align}
    &-(x+a)^2[\psi'(x+b)-\psi'(x+a)]+(a-b) \nonumber \\ 
    &=-(x+a)^2\left[
    \frac{a-b}{x^2}+\frac{(b-a)(a+b-1)}{x^3}+O\left(\frac{1}{x^4}\right)
    \right]+(a-b) \nonumber \\
    & =\frac{(a-b+1)(b-a)}{x}+O\left(\frac{1}{x^2}\right).\label{eq:psi_diff_asymp_finer}
\end{align}

Using \eqref{eq:psi_difference_asymp_finer} 
and \eqref{eq:psi_diff_asymp_finer} 
we obtain from \eqref{eq:upper_for_alpha_beta} 
\begin{align*}
    &\alpha(b-a-\beta)+\frac{\alpha(b-a)(a-b+1)}{2x}+O\left(\frac{1}{x^2}\right)\\
    &\leq b-a-\beta+\frac{(a-b+1)(b-a)}{x}+O\left(\frac{1}{x^2}\right).
\end{align*}
Letting $x\to\infty$ we conclude that $\alpha\leq 1$ when $\beta<b-a$, and 
$\alpha\leq 2$ when $\beta=b-a$.
 $\qedsymbol$

\textbf{Proof of Remark \ref{rem:alpha_less_2}.}

We know \cite{Luke1964} 2.11 (9)  that
\begin{align*}
	\psi(z+a) & =\log\,z- \sum_{k=0}^2
\frac{(-1)^{k+1}B_{k+1}(a)}{k+1}z^{-k-1}+O\left(\frac{1}{z^4}\right), \\
	& (z\to\infty\,\mathrm{in}\, |\arg\,z|<\pi-\varepsilon,\,\varepsilon>0),
\end{align*}
where $a$ is an arbitrary complex number, and $B_{k+1}(a)$ are Bernoulli polynomials, 
$B_1(a)=a-1/2$, $B_2(a)=a^2-a+1/6$, $B_3(a)=a^3-3a^2/2+a/2$  that is, 
\begin{equation*}
    \psi(x+a)=\log x+\frac{a-1/2}{x}-\frac{a^2-a+1/6}{2}\frac{1}{x^2}+
    \frac{a^3-3a^2/2+a/2}{3}\frac{1}{x^3}
    +O\left(\frac{1}{x^4}\right). 
\end{equation*}
Thus we get
\begin{align}
    & (x+a)[\psi(x+b)-\psi(x+a)] \nonumber\\
    & =(x+a)\left[ 
    \frac{b-a}{x}+\frac{(a-b)(a+b-1)}{2}\frac{1}{x^2} \right.\nonumber\\
    &\left.+\frac{(b-a)(1-3a+2a^2-3b+2ab+2b^2)}{6}\frac{1}{x^3}
    +O\left(\frac{1}{x^4}\right)\right] \nonumber\\
    &=b-a+\frac{(b-a)(a-b+1)}{2}\frac{1}{x}+\frac{(a-b)(1+a-b)(-1+a+2b)}{6}\frac{1}{x^2}
    \nonumber\\
    &+O\left(\frac{1}{x^3}\right).\label{eq:psi_difference_asymp_finest}
\end{align}
Since $\psi$ is analytical as $z\to\infty$ $(|\arg\,z|<\pi-\varepsilon,\,
\varepsilon>0)$, we obtain
\begin{align*}
	\psi'(z+a) & =\frac{1}{z}+ \sum_{k=0}^2
(-1)^{k+1}B_{k+1}(a)z^{-k-2}+O\left(\frac{1}{z^5}\right), \\
	& (z\to\infty\,\mathrm{in}\, |\arg\,z|<\pi-\varepsilon,\,\varepsilon>0),
\end{align*}
that is, 
\begin{equation*}
    \psi'(x+a)=\frac{1}{x}-\frac{a-1/2}{x^2}+\frac{a^2-a+1/6}{x^3}
    -\frac{a^3-3a^2/2+a/2}{x^4}
    +O\left(\frac{1}{x^5}\right). 
\end{equation*}
Thus we get
\begin{align}
    & (x+a)^2[\psi'(x+b)-\psi'(x+a)] \nonumber\\
    & =(x+a)^2\left[ 
    \frac{a-b}{x^2}+\frac{(b-a)(a+b-1)}{x^3} \right.\nonumber\\
    &\left. +\frac{(a-b)(1-3a+2a^2-3b+2ab+2b^2)}{2x^4}
    +O\left(\frac{1}{x^5}\right)\right] \nonumber\\
    & =a-b+\frac{(a-b)(1+a-b)}{x}+\frac{(b-a)(1+a-b)(-1+2b)}{2x^2}+O\left(\frac{1}
    {x^3}\right). \label{eq:psi_first_diff_finest}
\end{align}
Using \eqref{eq:psi_difference_asymp_finest} and \eqref{eq:psi_first_diff_finest} we obtain from 
\eqref{eq:upper_for_alpha_beta}
\begin{align*}
    &\alpha\{ (x+a)[\psi[x+b]-\psi(x+a)]-(b-a)\} \\
    &=\frac{\alpha(b-a)(a-b+1)}{2x}+\frac{\alpha(a-b)(1+a-b)(a-1+2b)}{6x^2}
    +O\left(\frac{1}{x^3}\right) \\
    &\leq -(x+a)^2[\psi'(x+b)-\psi'(x+a)]+(a-b) \\
    &=\frac{(b-a)(1+a-b)}{x}+\frac{(b-a)(1+a-b)(1-2b)}{2x^2}+O\left(\frac{1}{x^3}\right).
\end{align*}
Putting $2$ into $\alpha$ we obtain
\begin{equation*}
    (b-a)(1+a-b)(-1-2a+2b)\leq O\left(\frac{1}{x}\right),
\end{equation*}
which is impossible when $-1-2a+2b>0$, hence in this case $\alpha\neq 2$, 
that is $\alpha<2$. $\qedsymbol$

\textbf{Proof of Theorem \ref{thm:wrong_mon}.}
The function 
\begin{equation*}
	f_{0,\beta}:=(-\log\circ g_{\beta})'
\end{equation*}
is completely monotone by Theorem 
\ref{thm:psi_Chen2007_improve_comp_mon} when $\beta\leq b-a$.  

Now we prove the only if part. 
Assume that $g_{\beta}$ is logarithmically completely monotone. 
Corollary \ref{cor:log_f_diff_comp_mon} implies that $g_{\beta}$ is completely 
monotone.
Then for $x>-a$
\begin{align*}
	g_{\beta}'(x) &=g_{\beta}(x)(\log\circ g_{\beta})'(x)\\
	&=\frac{g_{\beta}(x)}{x+a}\left\{
	\beta+(x+a)\left[\psi(x+a)-\psi(x+b)\right]\right\}\leq 0,
\end{align*}
that is,
\begin{equation*}
	\beta\leq (x+a)\left[\psi(x+b)-\psi(x+a)\right].
\end{equation*}
Using that $\psi(x+b)-\psi(x+a)\geq 0$, the asymptotic formula 
\eqref{eq:psi_different} and applying $x\to\infty$ 
we obtain $\beta\leq b-a$.
 $\qedsymbol$

\textbf{Proof of Theorem \ref{thm:frac_gamma_power}.}
Define the function $\varphi$ by
\begin{equation*}
	\varphi(x):=\log\,\Gamma(x+b)-\log\,\Gamma(x+a)-\beta\log(x+a)
	+\log\,c_0.
\end{equation*}
Then 
\begin{equation*}
	\varphi'(x)=\psi(x+b)-\psi(x+a)-\frac{\beta}{x+a}=f_{0,\beta}(x).
\end{equation*}
Theorem \ref{thm:psi_Chen2007_improve_comp_mon} yields 
$\varphi'(x)$ is completely monotone on $[0,\infty)$.  Thus  
$\varphi'(x)\geq 0$ so $\varphi(x)$ is monotone increasing. 
Since $\varphi(0)\geq 0$ we obtain $\varphi(x)$ is a Bernstein function 
on $[0,\infty)$. Lastly, Theorem \ref{thm:f_powerto_g_log_comp} 
implies that our function is 
logarithmically completely monotone. $\qedsymbol$

\textbf{Proof of Theorem \ref{thm:log_gamma_vogt_general}.}
By Lemma \ref{lem:vogt_comp_mon} 
\begin{equation*}
	1+\frac{1}{x}\log\,\Gamma(x+1)-\log(x+1)
\end{equation*}
is completely monotone on $(-1,\infty)$.
By Theorem \ref{thm:linfraction_log_comp_mon}
\begin{equation*}
	\log\left(\frac{x+1}{x+\beta}\right)
\end{equation*}
is completely monotone on $(0,\infty)$. Thus their sum
\begin{equation*}
	1+\frac{1}{x}\log\,\Gamma(x+1)-\log(x+\beta)
\end{equation*}
is also completely monotone on $(0,\infty)$. $\qedsymbol$

\section{Acknowledgment}
The author is grateful to the referee for careful reading of the paper and for 
valuable comments, which helped to correct 
 errors and improve the quality of the paper.

\end{document}